\documentclass[12pt,leqno]{article}
\usepackage{amsfonts,amssymb,fullpage, epsfig}
\newtheorem{num}{\indent\hskip-4.5pt}[section]

\newcommand{\alku}{ \begin{num}\hskip-6pt \hskip5pt }
\newcommand{\loppu}{\end{num}}

\newcommand{\be}{\addtocounter{num}{1}\begin{equation}}
\newcommand{\ee}{\end{equation}}

\newcommand{\proof}{{\bf Proof. }\rm }

\begin{document}

\def\Rn{{\mathbb R}^n}
\def\Rnbar{\overline{{\mathbb R}^n}}
\def\Rb{{\overline \mathbb R}}
\def\R{{\mathbb R}}
\def\Rp{{\mathbb R}^+}
\def\C{{\mathbb C}}
\def\X{{\mathbb X}}
\def\cross{\times}
\def\ple{\preceq}
\def\pge{\succeq}
\def\ch{{\cosh}}
\def\sh{{\sinh}}
\def\arsh{{\rm arsh }}
\def\arch{{\rm arch }}
\def\lem{{\rm {\bf Lemma.\hskip 0.5truecm}}}
\def\pro{{\bf Proposition.\hskip 0.5truecm }}
\def\cor{{\bf Corollary.\hskip 0.5truecm}}
\def\thm{{\bf Theorem.\hskip 0.5truecm}}
\def\rem{{\bf Remark.\hskip 0.5truecm}}
\def\define{{\bf Definition.\hskip 0.5truecm}}
\def\ack{{\bf Acknowledgement.\hskip 0.5truecm}}
\def\lqq{\lq\lq}
\def\rqq{\rq\rq}
\def\rqqs{\rq\rq\hskip 0.15truecm}
\def\ineq{\not=}
\def\hbar{\vert}
\def\l{\ell}
\def\inv{ {\rm inv} } 
\def\abrv{.\hskip 0.1truecm}
\def\ident{\equiv}
\def\card{{\rm card\, }}
\def\GM{{GM(B^n)}}
% \sim?

\pagestyle{myheadings}
\title{A New Weighted Metric: the Relative Metric II}
\author{Peter A. H\"ast\"o\thanks{Supported in part by the Academy of Finland and the Finnish Academy of Sciences (Viljo, Yrj\"o and Kalle V\"aisal\"a's Fund)},
\\ 
{\small {Department of Mathematics, University of Helsinki, P.O. Box 4, 00014, Helsinki, Finland,}}
\\
{\small{E-mail: peter.hasto@helsinki.fi.}} }
\date{}
\maketitle

\abstract{\small{
In the first part of this investigation, \cite{Ha}, we generalized 
a weighted distance function of
\cite{Li} and found necessary and sufficient conditions for it being a metric. 
In this paper some properties of this so--called $M$--relative metric are 
established. Specifically, isometries, quasiconvexity and local convexity 
results are derived. We also illustrate connections between our approach and 
generalizations of the hyperbolic metric.}}

\bigskip

\centerline{Keywords: Relative metric, weighted metric}

\centerline{Mathematics Subject Classification (2000): Primary 39B62, Secondary 26D07, 30F45.} 

 %%%%%%%%%%%SECTION 1
 %%%%%%%%%%%SECTION 1
 %%%%%%%%%%%SECTION 1

\section{Preliminaries and main results}

In this section we introduce the $M$--relative metric and state the main results. In
order to do this, we have to introduce some notation -- for a fuller account 
the reader should consult Section 2 of \cite{Ha}.

A normed space $X$ is called {\it Ptolemaic} if 
$$ \Vert z-w\Vert \Vert x-y \Vert \le \Vert y -w\Vert \Vert x-z \Vert + \Vert x-w\Vert \Vert z-y \Vert$$
holds for every $x,y,z,w\in X$ (for background information on Ptolemy's inequality, see e\abrv g\abrv \cite[10.9.2]{Be}). Throughout this paper, we will denote by $\X$
a Ptolemaic normed space which is non--degenerate, i.e\abrv $\X$ is non--empty and 
$\X \ineq \{0\}$. By a metric or a norm we understand a function from $\X\cross \X$ 
into $[0, \infty]$ or $\X$ into $[0,\infty]$, respectively.

An increasing function $f \colon [0,\infty) \to [0,\infty)$ is said to be 
{\it moderately increasing} 
if $f(t)/t$ is decreasing. A function $P\colon [0,\infty)\cross [0,\infty)
 \to [0,\infty)$ of two variables 
is moderately increasing if both $P(x,\cdot)$ and $P(\cdot,x)$ are 
moderately increasing for each fixed $x\in [0,\infty)$. 

If $P\colon [0,\infty)\cross [0,\infty)\to [0,\infty)$ satisfies
$$\max\{x^{\alpha},y^{\alpha}\} \ge P(x,y) \ge \min\{x^{\alpha},y^{\alpha}\},$$
for all $x,y\in [0, \infty]$ then it is called an {\it $\alpha$--quasimean}. A 
1--quasimean is called a {\it mean}. We define the {\it trace} of a symmetric 
quasimean $P$ by $t_P(x):=P(x,1)$ for $x\in [1,\infty]$. We will need the following 
family of quasimeans 
$$ S_p(x,y):= (1-p){ x-y \over x^{1-p}-y^{1-p}},\ S_p(x,x)=x^p,\ 0<p<1, $$
$$ S_1(x,y):=L(x,y):={x-y\over \log x - \log y},\ S_1(x,x)=x.$$

Throughout this paper we will denote by $M$ a symmetric function, 
$M\colon [0,\infty]\cross [0,\infty) \to [0,\infty)$. When $M(x,y)=f(x)f(y)$ 
this means, then, that we assume that $f \colon [0,\infty) \to [0,\infty)$.
By the {\it $M$--relative distance} (in $\X$) we mean the function 
$$ \rho_M(x,y) := {\Vert x-y\Vert \over M(\Vert x \Vert, \Vert y\Vert)}$$
where $x,y\in \X$ (here we define 0/0=0). We will use the convention 
$ M(x,y):= M(\Vert x \Vert, \Vert y\Vert)$ (and $f(x):=f(\Vert x \Vert)$, when 
$M(x,y)=f(x)f(y)$). If $\rho_M$ is a metric, it is 
called the {\it $M$--relative metric}. The main results of the first part of 
this investigation are summarized in the next theorem.

\alku\thm\label{metricThm} {\rm (\cite[Sections 1 \& 3]{Ha})} 
Let $\X$ denote a non--degenerate Ptolemaic normed space.
\begin{itemize}
\item[{\rm (1)}] Assume that $M$ is moderately increasing. Then $\rho_M$ is a metric 
in $\X$ if and only if it is a metric in $\R$. 
\item[{\rm (2)}] Let $M$ is an $\alpha$--quasimean. Then $\rho_M$ is a metric in 
${\R}$ if $M(x,1)/S_{\alpha}(x,1)$ is increasing in $x$ for $x\ge 1$. If
$\rho_M$ is a metric in ${\R}$ then $M(x,1)\ge S_{\alpha}(x,1)$ for $x\ge 1$.
\item[{\rm (3)}] Assume that $M(x,y)=(x^p+y^p)^{q/p}$ for $p,q>0$. 
Then $\rho_M$ is denoted 
$\rho_{p,q}$ and called the ($p$, $q$)--relative distance. It is a metric 
in $\X$ if and only if $q=0$ or $0 < q \le 1$ and $p \ge \max\{ 1-q, (2-q)/3\}$.
\item[{\rm (4)}] Let $M(x,y)=f(x)f(y)$. Then $\rho_M$ 
is a finite metric (i.e\abrv $\rho_M<\infty$) 
in $\X$ if and only if $f$ is moderately increasing and convex.
\end{itemize}
\loppu

Like the first part of the investigation, this paper is organized along 
three threads -- one general and two special ones. 

In the general case, the moderation assumption also suffices for 
deriving some results on lipschitz mappings, quasiconvexity and local star--shapedness
of the metric (in Sections 2, 4 and 5, respectively). 

In the special cases, we can prove a bit more, however we also have to restrict 
ourselves to the spaces $\Rn$:

\alku\thm\label{pqMainThm} Let $\rho_{p,q}$ denote the 
{\it $(p,q)$--relative metric} as in Theorem \ref{metricThm} (3). Then 
\begin{itemize}
\item[{\rm (1)}] If $n \ge 2$, the $(p,q)$--relative metric is 
quasiconvex in $\Rn$ (see Section 4 for the definition) if and only if $q<1$ in which 
case it is $c_{p,q}$--quasiconvex, where
$$ {2^{-q/p} \over 1-q} \le c_{p,q} \le {\max\{ 2^{q(1-1/p)}, 1\} \over 1-q}.$$
\item[{\rm (2)}] The $(p,q)$--relative metric is locally convex 
(see Section 5 for the definition) if and only if $p<\infty$.
\end{itemize}
\loppu

\alku\thm\label{ffMainThm} Let $M(x,y)=f(x)f(y)$. If $n\ge 2$, $\rho_M$ 
is $c$--quasiconvex in $\X$ for some $c\le \sqrt{\pi^2/4+4}$.
\loppu

This paper also contains an explicit formula for the $\alpha$--quasihyperbolic 
metric in the domain $\Rn \setminus \{0\}$ which might be of independent interest 
(the $\alpha$--quasihyperbolic is defined in the begining of Section 4).

\alku\thm\label{k_alpha} For $n\ge 2$ and $0<\alpha< 1$ we have
$$ k_{\alpha}(x,y)= {1\over \beta} \sqrt{\vert x\vert^{2\beta} + \vert y\vert^{2\beta}
-2\vert x\vert^{\beta}\vert y\vert^{\beta}\cos \beta \theta}.$$
Here $\alpha+\beta =1$ and $\theta$ is the angle $\widehat{x0y}$.
In particular, as $\alpha\to 1$, 
$$ k_{\alpha}(x,y)\to \sqrt{\theta^2+ \log^2(\vert x\vert/\vert y\vert)},$$ 
the well--known expression for the quasihyperbolic
metric in $\Rn \setminus\{0\}$ {\rm (\cite[3.11]{Vu})}.
\loppu 

In the last section we consider how the relative-metric-approach may be applied 
to extending the hyperbolic metric in $\Rn$ for $n\ge 3$. We 
illustrate the limitations of the approach by considering a generalization of
the hyperbolic metric proposed in \cite[3.25, 3.26]{Vu}
concerning a metric similar to $\rho_M$ and proving the triangle inequality 
by another method.

%%%%%%%%%%%%%%%%%%%%%%%%%%%%%%%%%%%%%%%%%%%%%%%%%%%%%%%%%%%%%%%%%%%%%%%%%%%%%%%%%%%
%%%%%%%%%%%%%%%%%%%%%%%%%%%%%%%%%%%%%%%%%%%%%%%%%%%%%%%%%%%%%%%%%%%%%%%%%%%%%%%%%%%
%%%%%%%%%%%%%%%%%%%%%%%%%%%%%%%%%%%%%%%%%%%%%%%%%%%%%%%%%%%%%%%%%%%%%%%%%%%%%%%%%%%

\vskip 30pt
\section{Bilipschitz mappings and $\rho_M$}

\alku\lem\label{bilipLem} 
Let $M$ be moderately increasing, $\rho_M$ be a metric in $\X$ and
$g\colon \X \to \X$ be $L$--bilipschitz with respect to the norm $\Vert\cdot\Vert$ with 
$g(0)=0$. Then $g$ is $L^3$--bilipschitz with respect to the metric $\rho_M$.
\loppu

\proof Assume first that $x,y\ineq 0$. Since $M$ is increasing
$$ \rho_M(g(x),g(y))= {\Vert g(x)-g(y) \Vert \over M(g(x),g(y))} \le
{L \Vert x-y \Vert \over M(x/L,y/L)} \le L^3 {\Vert x-y \Vert \over M(x,y)},$$
where the last inequality follows since
$$ {M(x/L,y/L) \over xy/L^2} \ge {M(x,y/L) \over xy/L} \ge {M(x,y) \over xy },$$
by the moderation condition. On the other hand if $y=0$ and $M(g(x),0)>0$ then  
$$ \rho_M(g(x),0)= {\Vert g(x) \Vert \over M(g(x),0)} \le
{L \Vert x \Vert \over M(x/L,0)} \le L^2 {\Vert x \Vert \over M(x,0)}.$$
The case $M(g(x),0)=0$ is trivial and so the upper bound is proved. 
The lower lipschitz bound follows similarly. $\square$

\alku\rem{\rm It is clear that the condition $g(0)=0$ in Lemma \ref{bilipLem} is 
essential. For the translation $x\mapsto x+a$ is $1$--bilipschitz in the norm 
$\Vert\cdot\Vert$. If, for instance, $M(x,y)=x+y$ then 
$$ \lim_{\epsilon\to 0} {\rho_M(-\epsilon,\epsilon)\over \rho_M(a-\epsilon,a+\epsilon)} 
= \infty, $$
hence the translation is not bilipschitz in $\rho_M$. Note also that the 
condition $g(0)=0$ can be understood in terms of the generalization of the relative 
metric presented in Section 6 of \cite{Ha}: the $\rho_M$ is finite in 
$\X\setminus\{0\}$ if $M$ is moderately increasing (and $M\not\ident 0$) and hence 
the relevant class of mappings are from $\X\setminus\{0\}$ to $\X\setminus\{0\}$, 
i.e\abrv those with $g(0)=0$. 
}\loppu

\alku\lem\label{qcLem} Let $M$ be moderately increasing with $M\not\ident 0$ 
and let $\X=\Rn$. 
If $\rho_M$ is a metric and
$g\colon \Rn \to \Rn$ is $L$--bilipschitz with respect to the metric $\rho_M$ then
$g$ is quasiconformal in $\Rn$ with linear dilatation coefficient 
less than or equal to $L^2$.
\loppu

\proof We will first prove that $g$ is continuous in $\Rn\setminus\{0\}$. 
Since $M$ is moderately continuous 
and $M\not\ident 0$ it follows that $M(x,y)>0$ unless $xy=0$. 

Fix a point $x\in\Rn$ such that $g(x)\ineq 0$. Since $\rho_M$ is a metric and $g$ is 
bilipschitz with respect to $|rho_M$ it follows that $g$ is injective. Hence there exists 
a neighborhood $U_x$ of $x$ such that $0\not\in U_x$ and $g(y)\ineq 0$ for $y\in U_x$.    
For $y\in U$ $\Vert g(y)\Vert$ has an upper bound independent of $y$. 
For if $\Vert g(y)\Vert \ge \Vert g(x)\Vert $ then the inequality
$$ \Vert g(y)\Vert -\Vert g(x)\Vert  \le \Vert g(x)-g(y)\Vert  \le L{\Vert g(y)\Vert \over \Vert g(x)\Vert }M(g(x),g(x)) \rho_M(x,y)$$
implies that 
$$ \Vert g(y)\Vert  \le \Vert g(x)\Vert \left(1- {L\over \Vert g(x)\Vert }M(g(x),g(x))\rho_M(x,y)\right)^{-1}. $$
It follows that  
$$ \Vert g(x)-g(y)\Vert  \le {L M(g(x),g(x))\over 1- LM(g(x),g(x)) \rho_M(x,y)/\Vert g(x)\Vert }
\rho_M(x,y). $$
From this we easily see that $g(y)\to g(x)$ as $y\to x$. Hence $g$ is 
continuous in $\Rn\setminus\{g^{-1}(0)\}$. 

Let $z\not\in \{0,g^{-1}(0)\}$, $x,y\in U_z$ and 
$\Vert x-z\Vert = \Vert y-z\Vert =r$. Then
$$ {\Vert g(x)-g(z) \Vert \over \Vert g(y)-g(z) \Vert} \le L^2{M(g(x),g(z))M(y,z)\over
M(g(y),g(z))M(x,z)}.$$
Since $M$ is moderately increasing it is continuous in $\Rn\setminus\{0\}$ by 
\cite[Lemma 2.3]{Ha}. By the continuity of $M$ and $g$ the right hand side tends 
to $L^2$ as $r\to 0$. Hence we have proved that $g$ is quasiconformal in 
$\Rn\setminus\{0,g^{-1}(0)\}$. But then $g$ is quasiconformal in $\Rn$ by 
well-known continuation results (see e.g\abrv \cite{Va2}). 
$\square$.

\alku\rem{\rm If $M$ and $g$ are as in the previous lemma and additionally 
$M(x,0)=0$ for every $x>0$ then $g(0)=0$. 
For the bilipschitz condition 
$$ {1\over L} {\Vert x-y\Vert \over M(x,y)} \le {\Vert g(x)-g(y)\Vert  \over M(g(x),g(y))} \le 
L {\Vert x-y\Vert \over M(x,y)}$$
implies that $M(x,y)$ and $M(g(x),g(y))$ are simultaneously $0$. 
Therefore $\Vert x\Vert \Vert y\Vert =0$ iff $M(x,y)=0$ iff $M(g(x),g(y))=0$ iff 
$\Vert g(x)\Vert \Vert g(y)\Vert =0$, which implies $g(0)=0$.
}\loppu 

\alku\cor\label{confCor} If $M$ is moderately increasing with $M\not\ident 0$ 
and $g\colon \X \to \X$ is a $\rho_M$--isometry then $g$ is conformal. $\square$
\loppu

\alku\rem {\rm The mapping $g(x)=\vert x\vert x$ is 2--bilipschitz in the 
$\rho_{\infty,1}$ metric (=$\rho_M$ with $M(x,y)=\max\{x,y\}$) 
but is not lipschitz with respect to the Euclidean metric 
(=$\rho_M$ with $M\ident 1$). The spherical metric, $q$ (=$\rho_M$ with 
$M(x,y)=\sqrt{1+x^2}\sqrt{1+y^2}$) and the inversion 
$x \mapsto x/\Vert x\Vert^2$ is a $q$--isometry. However, this inversion is 
certainly not lipschitz with respect to the Euclidean metric. These examples 
show that the class of $\rho_M$--lipschitz mappings depends on $M$ in a non--trivial way.
}\loppu

%%%%%%%%%%%%%%%%%%%%%%%%%%%%%%%%%%%%%%%%%%%%%%%%%%%%%%%%%%%%%%%%%%%%%%%%%%%%%%%%%%%
%%%%%%%%%%%%%%%%%%%%%%%%%%%%%%%%%%%%%%%%%%%%%%%%%%%%%%%%%%%%%%%%%%%%%%%%%%%%%%%%%%%
%%%%%%%%%%%%%%%%%%%%%%%%%%%%%%%%%%%%%%%%%%%%%%%%%%%%%%%%%%%%%%%%%%%%%%%%%%%%%%%%%%%

\vskip 30pt 
\section{$\alpha$--quasihyperbolic metrics}

The length of a (rectifiable) path $\gamma \colon [0,l] \to \X$ in the 
metric $\rho_M$ with continuous $M$ is defined by
$$ \l_M(\gamma) := \lim_{n\to\infty} 
\sum_{i=0}^{n} \rho_M(\gamma(t_i), \gamma(t_{i+1})),$$
where $t_i<t_{i+1}$, $t_0=0$, $t_n=l$ and $\max \{t_{i+1}-t_i\} \to 0$. 
If $\gamma$ is any path connecting $x$ and $y$ in $\X$ then 
$\rho_M(x,y)\le \l_M(\gamma)$ by the triangle inequality. 

Let $M$ be an $\alpha$--quasimean ($0<\alpha\le 1$). By taking the infimum over 
all rectifiable paths joining $x$ and $y$ we conclude that
$$ \rho_M(x,y)\le \inf_{\gamma} \l_M(\gamma) = \inf_{\gamma} \int {ds\over 
\Vert \gamma(s)\Vert^{\alpha}}
=:k_{\alpha}(x,y),$$
since $M(x,x+\epsilon) \ge x^{\alpha}$ for $\epsilon>0$.
Here $k_\alpha$ stands for the $\alpha$--quasihyperbolic metric, which was 
introduced in \cite{GM}. More precisely, it is the $\alpha$--quasihyperbolic metric in
the domain $G=\Rn \setminus \{0\}$. In this section we will 
derive an explicit expression for $k_{\alpha}(x,y)$, which will be used to study 
quasiconvexity in the next section. 

%\alku\thm\label{k_alpha} For $n\ge 2$, $0<\alpha\le 1$ and 
%$\vert x \vert \ge \vert y\vert$ we have
%$$ k_{\alpha}(x,y)= {1\over \beta} \sqrt{\vert x\vert^{2\beta} + \vert y\vert^{2\beta}
%-2\vert x\vert^{\beta}\vert y\vert^{\beta}\cos \beta \theta}.$$
%Here $\alpha+\beta =1$ and $\theta$ is the angle $\widehat{x0y}$.
%In particular, as $\alpha\to 1$, $k_{\alpha}(x,y)\to \sqrt{\theta^2+
%\log^2(\vert x\vert/\vert y\vert)}$, the well--known expression for the quasihyperbolic
%metric in $\Rn \setminus\{0\}$ {\rm (\cite[3.11]{Vu})}.
%\loppu 

% replaced reference
\alku{\bf Proof of Theorem 1.4. \hskip 0.5truecm } {\rm It is clearly sufficient to limit ourselves to the
case $\X=\R^2$ in this proof. It is also clear that the geodesic can be parameterized by 
$(r(\theta),\theta)$ in polar coordinates. The kernel of the integral then becomes
$r^{-\alpha}\sqrt{(r')^2+r^2}$, where $r'=dr/d\theta$. Then the Euler equation (cf\abrv \cite[p\abrv 36 (5)]{El}) tells 
us that the geodesic satisfies the differential equation
$$ r^{-\alpha}\sqrt{(r')^2+r^2}- { r^{-\alpha} (r')^2 \over \sqrt{(r')^2+r^2} } =c_1.$$
Since $c_1$ is independent of $r$, one easily sees that $c_1\ineq 0$. 
Then the equation is equivalent to $r^{\beta}/c_1=\sqrt{((\log r)')^2+1}$. 

To solve this equation, we change variables by substituting $y:=\log r$. 
The equation then becomes $e^{\beta y} = c_1\sqrt{(y')^2+1}$, where $y'=dy/d\theta$. 
We introduce an auxiliary parameter, $t$, by $\sinh t =y'$. 
Then $ e^{\beta y}=c_1\cosh t$ and 
$$ d\theta = {dy/dt \over dy/d\theta}dt= { dt\over \beta \cosh t}.$$

Solving this equation gives $\tan ((\beta \theta + c_2)/2) = e^t$, hence 
$$ r(\theta)^{\beta}= {c_1\over 2} \left(\tan ((\beta\theta + c_2)/ 2) + 
{1 \over \tan ((\beta \theta + c_2)/2) } \right) = 
{c_1\over \sin (\beta \theta + c_2)}.$$

Let us now calculate the distance in the $k_\alpha$ metric between $1$ and $re^{i\theta_1}$, where $r\ge 1$ and $0\le \theta_1\le \pi$, using the formula for the geodesic (denoted by $\gamma$):
$$ k_\alpha(1, re^{i\theta_1}) = \int_{\gamma} {\sqrt{(r')^2+r^2}\over r^{\alpha}}d\theta
= \int_0^{\theta_1} {c_1\over \sin^2(\beta \theta + c_2)} d\theta = 
{c_1\over \beta} \left(\cot c_2-\cot(\beta \theta_1 + c_2) \right).$$
It remains to express $c_1$ and $c_2$ in terms of the boundary values: 
$$ \sin c_2 = c_1,\ r^{\beta}\sin ( \beta\theta_1+c_2)=c_1.$$
These equations imply that 
$$ c_1= {r^\beta \sin \beta\theta_1 \over \sqrt{ 1+r^{2\beta} - 
2 r^{\beta}\cos\beta\theta_1}},$$
from which it follows that
$$ k_\alpha(1, re^{i\theta_1}) = \frac{1}{\beta}\left(\sqrt{r^{2\beta}-c_1^2} 
\pm \sqrt{1-c_1^2}\right) = {r^\beta \vert r^\beta - \cos \beta\theta_1 \vert 
\pm \vert r^\beta \cos \beta\theta_1 - 1 \vert \over \beta \sqrt{1+r^{2\beta} - 
2r^{\beta}\cos\beta\theta_1}},$$
where $\pm$ is a plus when $c_2$ is greater than $\pi/2$ and a minus when it is
not. This means that effectively the absolute value is disregarded and the $\pm$ sign
is a minus sign since $c_2$ is greater than $\pi/2$ exactly when 
$ r^\beta \cos \beta\theta_1 \ge 1$. 

Then 
$$ r^\beta \vert r^\beta - \cos \beta\theta_1 \vert \pm \vert r^\beta 
\cos \beta\theta_1 - 1 \vert = r^\beta (r^\beta - \cos \beta\theta_1) 
- (r^\beta \cos \beta\theta_1 - 1) = $$
$$ = 1+r^{2\beta} - 
2r^{\beta}\cos\beta\theta_1$$
from which the claim follows. $\square$
}\loppu

\alku\rem{\rm To get a picture of what $k_\alpha$ looks like we consider how
the distance between points changes as $\alpha$ changes. Since $k_\alpha$ is
$\beta$-homogeneous and spherically symmetric, we assume that $y=1$. Consider
first the case when $x$ is a real number greater than one. Then $k_\alpha(x,1) = 
(x^{\beta}-1)/\beta$. This is an increasing function with respect to $\beta$.
Consider now another point $z\in S^{n-1}(0,1)$. Then 
$k_\alpha(z,1)= \sqrt{2(1-\cos\beta\theta)}/\beta$. This is decreasing 
in $\beta$. Hence, intuitively speaking, increasing $\alpha$ increases angular 
distance but decreases radial distance. Note that these considerations imply, 
in particular, that $k_\alpha$ is not monotone in $\alpha$.}\loppu

\alku\cor\label{KaEst} Let $\alpha + \beta =1$ with $0\le\alpha<1$. Then
$ k_\alpha(x,y)\le (\vert x \vert^\beta + \vert y \vert^\beta)/\beta.\ \square$
\loppu

\alku\lem\label{decKa} Let $\alpha + \beta =1$ with $0\le\alpha<1$. Then
$$ {\vert x\vert^\beta - \vert y\vert^\beta \over 
\beta(\vert x\vert - \vert y\vert)} \le {k_{\alpha}(x,y)\over \vert x-y\vert} \le 
{k_{\alpha}(-\vert x\vert,\vert y\vert) \over\vert x\vert+\vert y\vert} 
\le {\vert x\vert^\beta + \vert y\vert^\beta \over \beta(\vert x\vert + \vert y\vert)} 
\le {2^{\alpha} \over \beta} \vert x-y\vert^{-\alpha}.$$
Let $t\in\Rp$. There is equality in the first inequality for $x=ty$, in the second 
for $x=-ty$, in the third for $x=-ty$ and $\beta=1$ and in the fourth for $x=-y$.  
\loppu

\proof It suffices to show that 
$$ {k_{\alpha}(re^{i\theta},1) \over \vert re^{i\theta}-1 \vert }$$
is increasing in $\theta$ for $r\ge 1$ and $0\le \theta\le \pi$. 
Using the explicit formula for $k_\alpha$ from 
Theorem \ref{k_alpha} we need to show that
$$ {1+r^{2\beta}-2r^{\beta}\cos \beta\theta\over \beta(1+r^2 -2r\cos\theta)}$$
is increasing in $\theta$. We differentiate the equation with respect to $\theta$ and 
see that this follows if we show that
$$(1+r^{2\beta}-2r^{\beta}\cos \beta\theta)/(\beta r^{\beta}\sin\beta\theta)$$
is increasing in $\beta$. 

When we differentiate this equation with respect to $\beta$, we see that it suffices 
to show that
\be\label{decKa1} (s - 1/s)\log s\sin x + 2x + \sin 2x \ge (s+1/s)(x\cos x + \sin x),\ee
where we have denoted $s:= r^{\beta}\ge 1$ and $0\le x:=\beta\theta \le \pi$.
The inequality holds in (\ref{decKa1}) for $s=1$ since $x-\sin x \ge 
\cos x (x-\sin x)$ for $x\ge 0$. Differentiating (\ref{decKa1}) with respect to 
$s$ leads to
$$ (s^2+1)\log s +s^2-1 \ge (s^2-1)(x/\tan x + 1).$$
Since $x/\tan x\le 1$ for $0\le x\le \pi$, it suffices to show that $\log s \ge 1 - 2/(s^2+1)$, which follows 
since $2/(s^2+1)+\log s$ is increasing in $s$. $\square$ 

\alku\rem{\rm It would be interesting to see how the above estimates for $k_\alpha$ 
generalize to other domains than $\Rn\setminus \{0\}$.}\loppu

%%%%%%%%%%%%%%%%%%%%%%%%%%%%%%%%%%%%%%%%%%%%%%%%%%%%%%%%%%%%%%%%%%%%%%%%%%%%%%%%%%%
%%%%%%%%%%%%%%%%%%%%%%%%%%%%%%%%%%%%%%%%%%%%%%%%%%%%%%%%%%%%%%%%%%%%%%%%%%%%%%%%%%%
%%%%%%%%%%%%%%%%%%%%%%%%%%%%%%%%%%%%%%%%%%%%%%%%%%%%%%%%%%%%%%%%%%%%%%%%%%%%%%%%%%%

\vskip 30pt 
\section{Quasiconvexity}

In this section, we will assume that $n\ge 2$ and consider the space $\Rn$. The length
of a curve was defined at the beginning of the previous section. Following \cite{Va}, 
we define a metric $\rho_M$ (actually a metric space, ($\X$, $\rho_M$)) to 
be $c$--quasiconvex if $\inf_{\gamma} \l_M(\gamma) 
\le c \rho_M(x,y)$, where that infimum is taken over all rectifiable paths $\gamma$ 
joining $x$ and $y$. For instance if $G\subset \Rn$ is convex then ($G$, $|\cdot |$) is 
$1$--quasiconvex, whereas ($D$, $|\cdot |$) is not quasiconvex for $D:=B^n \setminus [0,1)$, since we need a path of lenght $\ge 1$ to connect $x:= (1-t)e_1 + te_2$ with 
$x:= (1-t)e_1 - te_2$ ($e_1$ and $e_2$ are basis vectors of $\Rn$).

\alku\thm\label{qconvM0} Let $M$ be an $\alpha$--quasimean such that
$\rho_M$ is a metric. Then $\rho_M$ is quasiconvex if and only if
$$ c_M := \sup_{x\ge 0,y>0}{k_\alpha(x,-y) \over x+y}M(x,y) <\infty, $$
in which case it is $c_M$--quasiconvex.
\loppu

\proof The claim follows directly from the second inequality in Lemma \ref{decKa}, since 
$\inf_{\gamma} \l_M(\gamma) = k_\alpha(x,y)$ by definition. $\square$.

\alku\cor\label{qconvMCor1} Let $M$ be $\alpha$--homogeneous with $M(1,1)=1$ such that
$\rho_M$ is a metric. Then $\rho_M$ is quasiconvex if and only if
$$ c_M := \sup_{r\ge 1}{k_\alpha(r,-1) \over r+1} M(r,1) <\infty, $$
in which case it is $c_M$--quasiconvex. $\square$
\loppu 

\alku\cor\label{qconvMCor2} Let $M$ be $\alpha$--homogeneous with $M(1,1)=1$ such that
$\rho_M$ is a metric. If $\alpha<1$ then $\rho_M$ is 
$2^\alpha /(1-\alpha)$--quasiconvex. If $M\le A_p^{\alpha}$ then $\rho_M$ is $c_{p,\alpha}$--quasiconvex, where
$$ c_{p,\alpha}:= {\max\{ 2^{\alpha(1-1/p)}, 1\} \over 1-\alpha}.$$
\loppu

\proof Let us consider $M=A_1^{\alpha}$. Then, by Corollary \ref{qconvMCor1} 
and Lemma \ref{decKa},
$$ c_M \le \sup_{r\ge 1}{ r^{1-\alpha}+1 \over (1-\alpha)(r+1)}
\left({r+1\over 2}\right)^{\alpha}= 
{1\over 2^{\alpha}(1-\alpha)} \sup_{r\ge 1}{r^{1-\alpha}+1 \over (r+1)^{1-\alpha}} 
\le 1/(1-\alpha),$$ 
since $(r^{1-\alpha}+1)(r+1)^{\alpha-1}$ is decreasing. 

Since $A_p\le \max\{ 2^{1-1/p},1\} A_1$ the second claim is proved. Since 
$$ M(x,1)\le A_{\infty}^{\alpha}(x,1)\le \{2A_1(x,1)\}^{\alpha} $$
for every $\alpha$--homogeneous $M$, the first claim also follows. $\square$.

%\alku\cor\label{pqQconv} The $(p,q)$--relative metric is quasiconvex if and only
%if $q<1$ in which case it is $c_{p,q}$--quasiconvex, where
%$$ {2^{-q/p} \over 1-q} \le c_{p,q} \le {\max\{ 2^{q(1-1/p)}, 1\} \over 1-q}.$$
%\loppu

% replaced reference 
\alku{\bf Proof of Theorem 1.2(1). \hskip 0.5truecm } {\rm 
The upper bound follows from Corollary \ref{qconvMCor2}. For the lower
bound let $r\to \infty$ in Corollary \ref{qconvMCor1}. $\square$
}\loppu

\alku\cor\label{qconv_q=1/2} $\rho_{p,1/2}$ is $\max\{\sqrt{2}, 2^{1-1/(2p)}\}$--quasiconvex, where the constant is the smallest possible.
\loppu

\proof Setting $\alpha=1/2$ in Corollary \ref{qconvMCor1} yields 
$$c_M= \sup_{r\ge 1} 2^{1-1/(2p)}\sqrt{(r^p+1)^{1/p}/(r+1)}, $$
from which the claim follows since $(r^p+1)^{1/p}/(r+1)$ is increasing for 
$p\ge 1$ and decreasing for $p\le 1$. $\square$

%\alku\thm\label{qconvff} Let $M(x,y)=f(x)f(y)$ be such that $\rho_M$ is a metric.
%Then $\rho_M$ is $c$--quasiconvex, where $c\le \sqrt{\pi^2/4+4}$.
%\loppu

% replaced reference
\alku{\bf Proof of Theorem 1.3. \hskip 0.5truecm } {\rm 
We will handle the cases $f(0)=0$ and $f(0)>0$ separately. In the first case 
$f(x)=cx$ for some $c$, as was shown in \cite[Remark 5.1]{Ha}. Denote by $x'$ 
the image of $x$ under the inversion $x\mapsto x/|x|^2$. Then 
$\rho_M(x,y)=|x'-y'|$ and hence the line from $x'$ to $y'$ is mapped onto a curve 
$\gamma$ (actually a segment or an arc of a circle) with $\l_M(\gamma) = \rho_M(x,y)$, 
hence $\rho_M$ is $1$--quasiconvex in this case. 

In the second case we may assume without loss of generality that $f(0)=1$. 

Let us fix the points $x$ and $y$ with $\Vert x\Vert\ge \Vert y\Vert>0$. Denote
by $\gamma_1$ the path which is radial from $x$ to $ (\Vert y\Vert /\Vert x\Vert ) x$
and then circular (with radius $\Vert y\Vert$) about the origin to $y$ and 
by $\gamma_2$ the path which is first circular (with radius $\Vert x\Vert$) and then radial from $(\Vert x\Vert /\Vert y\Vert ) y$ to $y$. 

In what follows we will denote $\Vert x \Vert$ by $x$ and similarly for $y$ and $z$, since there is no danger of confusion. 
We derive estimates for the lengths of the $\gamma_i$:
$$\min\{\l_M(\gamma_1), \l_M(\gamma_2)) \le \theta\min \left( {x 
\over f(x)^2}, {y \over f(y)^2}\right) + 
\int_{\Vert y \Vert}^{\Vert x \Vert } {dz\over f(z)^2}, $$
where $\theta$ is the angle $\widehat{x0y}$. 
Since $f$ is moderately increasing and convex we find that 
$f(z)\ge \max \{ 1 + z (f(y)-1)/y , z f(x)/x \}$ for $z\in[y , x ]$. 
Let $z_0\in Rp$ be such that $1 + z_0 (f(y)-1)/y = z_0 f(x)/x$. Then 
$$\int_{\Vert y \Vert}^{\Vert x \Vert } {dz\over f(z)^2} \le 
\int_{\Vert y \Vert}^{z_0} {dz\over \{1 + z(f(y)-1)/y  \}^2} + 
\int_{z_0}^{\Vert x \Vert} {dz\over \{z f(x)/x \}^2}\le $$
$$ \le {2\Vert x \Vert\over f(x)} - {y\over f(y)}- {\Vert x \Vert\over f(x)^2}
\left({x\over y}(f(y)-1)+1\right)\le {2(\Vert x \Vert-y)\over f(x)f(y)}.$$ 
To see that the last inequality holds, multiply by $f(x)^2f(y)$ and rearrange:
$$ 2(\Vert x \Vert-y)f(x) - 2\Vert x \Vert f(y)f(x) + yf(x)^2 \ge \left({\Vert x \Vert \over y}(f(y)-1)+1\right) \Vert x \Vert  f(y).$$
Notice that the right hand side is independent of $f(x)$ whereas the 
left hand side is increasing in $f(x)$ since 
$$y(f(x)-1)= (y-0)(f(x)-f(0))\ge (\Vert x \Vert -0)(f(y)-f(0))= \Vert x \Vert (f(x)-1),$$
which follows from the convexity of $f$. The inequality then follows, when 
we insert the minimum value for $f(x)$, that is $\Vert x \Vert (f(y)-1)/y+1$ and use 
$y(f(x)-1)\ge \Vert x \Vert (f(y)-1)$ again. 

In the case $y=0$ which was excluded above one easily derives the estimate 
$$ \l_M(\gamma_1) \le {2f(x)-1\over f(x)^2} x \le {2x\over f(x)}.$$

Now $c$--quasiconvexity follows, if we show that
$$ \theta \min\left(x{f(y)\over f(x)}, y {f(x)\over f(y)} \right) + 2(x-y)\le 
c\sqrt{x^2+y^2-2xy\cos\theta}.$$
For fixed $x$ and $y$, $\min\{x f(y)/f(x),y f(x)/f(y)\}\le \sqrt{xy}$. Hence it suffices 
to show that 
$$ \theta^2xy+4\theta(x-y)\sqrt{xy}+4(x-y)^2+2c^2xy\cos\theta \le c^2(x^2+y^2).$$
Since the case $y=0$ is clear we set $s:=x/y\ge 1$ and divide through by $xy$, obtaining:
$$ \theta^2+4(\sqrt{s}-\sqrt{1/s})\theta+4(\sqrt{s}-\sqrt{1/s})^2+2c^2\cos\theta - c^2(s+1/s) \le 0.$$
The derivative of the left  hand side with respect to $s$ is positive when
$$ 2\theta (s+1) \ge (c^2-4)\sqrt{s}(s-1/s)$$
or, equivalently, when $ \sqrt{s}-\sqrt{1/s}\le 8\theta/\pi^2$. Hence the only zero
of the derivative is a maximum, and we have 
$$ \theta^2+4(\sqrt{s}-\sqrt{1/s})\theta+4(\sqrt{s}-\sqrt{1/s})^2+2c^2\cos\theta - c^2(s+1/s) \le $$
$$ \le (1+16\pi^{-2})^2 \theta^2 +2c^2\cos\theta - 2c^2(32\theta^2/\pi^4 +1).$$
To see that the last expression in the inequality is less than zero, 
we use the expression $\pi^2/4+4$ for $c^2$:
$$ (1+16\pi^{-2})^2\theta^2 +2 (\pi^2/4+4)(\cos\theta -32\theta^2 \pi^4-1) \le 0.$$
When we divide by $1+16\pi^{-2}$, we see that this is equivalent to $\theta^2\le \pi^2(1-\cos\theta)/2$, which concludes the proof. $\square$
}\loppu

\alku\rem{\rm The first part of the proof of the previous theorem shows that for 
the universal constant $c$ for which every $\rho_M$ with $M(x,y)=f(x)f(y)$ is 
$c$--quasiconvex is at least $2$. For if $x$ and $y$ are on the same ray emanating 
from the origin the clearly the segment of the ray between $x$ and $y$ is the geodesic. 
Moreover the above derivation up to 
$$\int_{\Vert y \Vert}^{\Vert x \Vert } {dz\over f(z)^2} \le {2(x-y)\over f(x)f(y)}$$
is sharp. Hence $c\ge 2$, as claimed.
}\loppu
 
Metrics that are 1--quasiconvex are particularly interesting, since in these metric 
spaces any two points can be connected with a a path $\gamma$ 
with $\l_M(\gamma) = d(x,y)$, where $d$ is the metric, which is to say that the 
metric equals its own inner metric. The next 
lemma shows that, except for the Euclidean distance and its \lqq reciprocal\rqq, 
there are no 1--quasiconvex $M$--relative metrics in $\Rn$ with $n\ge 2$.

\alku\lem Let $M$ be moderately increasing. Then $\rho_M$ is a 1--quasiconvex 
metric in $\Rnbar$ if and only if $M\ident c>0$ or $M(x,y)=cxy$. 
\loppu

\proof In this proof we will write $r$ for $re_1$ etc. If $M\ident c>0$ or $M(x,y)=xy$ 
then clearly $\rho_M$ is 1--quasiconvex (the latter claim was shown in the 
Proof of Theorem 1.3). Assume conversely that $\rho_M$ is 1--quasiconvex. 
Consider the $1$--quasiconvex path $\gamma$, connecting $-r$ and $r$, where $r>0$. 

Now either $\infty\in\gamma$ or $\gamma$ crosses the $e_2$-axis. In the latter 
case let $b\in [0,\infty)$ be such that $\gamma$ crosses the $e_2$-axis in $b e_2$. Then, 
by the triangle (in)equality,  
$$ {2r\over M(r,r)} = {2\sqrt{r^2+b^2} \over M(r,b)} $$
or, equivalently, $M(r,b)=\sqrt{1+(b/r)^2} M(r,r)$. Suppose that $b\ineq 0$. 
Then $M(r,b)>M(r,r)$ and $b>r$ since $M$ is increasing and hence $(b/r)M(r,r)\ge M(r,b)$
since $M$ is moderately increasing. It follows that  
$$ {b\over r}M(r,r)\ge M(r,b) = \sqrt{1+(b/r)^2} M(r,r)$$
from which it follows that $b/r=\sqrt{1+(b/r)^2}$, which is impossible, hence $b=0$. 

It then follows that the path connecting $-r$ and 
$r$ is the segment $[-r,r]$. By considering the triangle equality for 
a point $a$, with $a<r$, on the path we find that $M(r,a)=M(r,r)$. We then consider 
again three distinct points $y$, $z$ and $x$ on $[0,r)$ in this order. 
The triangle equality becomes 
$$ {|x-y|\over M(x,x)}= {|x-z|\over M(x,x)} + {|z-y|\over M(z,z)},$$
hence $M(x,x)=M(z,z)$. But then $M(x,y)=M(x,x)=M(z,z)=M(z,w)$ (assuming $x\ge y$ and 

$z\ge w$, similarly otherwise) and we conclude $M(x,y)=c$ for $x,y\le r$.

Hence $\gamma$ does not cross $e_2$-axis, and we have $\infty\in\gamma$. 
This means that the path is the segment $[-\infty, -r]\cup [r,\infty]$. Now we may 
choose any point $b$, with $b\ge r$ on the path and get $2r/M(r,r)=2b/M(r,b)$, hence
$M(r,b)=(b/r)M(r,r)$ for all $b\ge r$. Then consider three arbitrary distinct 
points $y$, $z$ and $x$ on $(r,\infty)$ in this order. The triangle equality becomes
$$ {y-y^2/x \over M(y,y)} = {x-y\over M(x,y)} = {x-z\over M(x,z)}+ {z-y\over M(z,y)} =
{z-z^2/x\over M(z,z)} + {y-y^2/z\over M(y,y)}.$$
This leads to 
$$ {y^2/z-y^2/x\over M(y,y)} = {z^2/z - z^2/x\over M(z,z)},$$
hence, since $1/z-1/y\ineq 0$, $M(y,y)/y^2=M(z,z)/z^2$ for $y<z$. It then follows that
$M(r,b)=(b/r)M(r,r) = br M(1,1)$, i.e\abrv $M$ is of the form $M(x,y)=cxy$ for all $x,y\ge r$.

We have seen that there are two possible cases, either $M(x,y)=c$ for every 
$x,y\in B^2(0,r)$ or $M(x,y)=cxy$ for every $x,y\not\in B^2(0,r)$. If there is a 
path from $-r$ to $r$ trough $0$ then the same path will connect $-r'$ with 
$r'$ for $r'<r$ as well. Similarly for paths through $\infty$ and $r'>r$. Hence 
there exists an $r_0$ such that $M(x,y)=c$ for $|x|,|y|\le r$ and 
$M(x,y)= cxy/r^2$ for $|x|,|y|\ge r$. If $r_0=0$ or $r_0=\infty$ then everything 
$M$ equals $cxy$ or $c$ in the whole space.

Assume then that $0<r_0<\infty$. We may assume without loss of generality that 
$c=r_0=1$. Consider then the points $1/2$ and $2$. The $1$--quasiconvex path connecting 
these points goes through $1$, hence 
$$ \rho_M(1/2,2)= \rho_M(1/2,1) + \rho_M(1,2) = 1/2 + 1/2=1$$
and $ M(1/2,2) = (3/2)\rho_M(1/2,2)=3/2$. The $1$--quasiconvex path connecting 
$-1/2$ with $2$ crosses $S^{n-1}(0,1)$ at some point $z$. 
If $\theta=\widehat{2 0 z}$ then 
$$ \rho_M(2,z)= \sqrt{5/4 - \cos\theta},\ \rho_M(-1/2,z)= \sqrt{5/4 - \cos\theta},$$
so that $\rho_M(2,z) + \rho_M(-1/2,z) \ge 2> 5/3 = (5/2)/ M(1/2,2) = \rho_M(-1/2,2),$
contrary to the assumtion that $z$ lies on a $1$--quasiconvex path. This contradiction 
shows that this mixed case cannot occur. $\square$

\alku\rem{\rm Note that the question of when a generalized relative metric, of the 
type introduced in Section 6 of \cite{Ha} are quasiconvex is not directly answered 
by the results in this section. However since the quasiconvexity of either the 
$j_G$ metric or Seittenranta's metric, which are both generalized relative 
metrics, characterize uniform domains this question is clearly of interest. 
(See \cite[4.3-4.5]{Se}.)}
\loppu

%%%%%%%%%%%%%%%%%%%%%%%%%%%%%%%%%%%%%%%%%%%%%%%%%%%%%%%%%%%%%%%%%%%%%%%%%%%%%%%%%%%
%%%%%%%%%%%%%%%%%%%%%%%%%%%%%%%%%%%%%%%%%%%%%%%%%%%%%%%%%%%%%%%%%%%%%%%%%%%%%%%%%%%
%%%%%%%%%%%%%%%%%%%%%%%%%%%%%%%%%%%%%%%%%%%%%%%%%%%%%%%%%%%%%%%%%%%%%%%%%%%%%%%%%%%

\vskip 30pt 
\section{Local convexity}

In this section we consider how the relative metric grows in different directions. We
will denote by $B_d(x,r):= \{ y\in \X \colon d(x, y) <r\}$ denote the open ball in 
the metric space $(\X,d)$ and by $B^n(x,r)$ the Euclidean open ball of radius 
$r$ centered at $x$. Also $S_d(x, r)=\partial B_d(x, r)$ and 
$S^{n-1}(x, r)=\partial B^n(x, r)$. We will use the abbreviation $B_{\rho_M} =: B_M$ and
$S_{\rho_M} =: S_M$.

\alku\define {\rm 
\begin{itemize}
\item[(i)] 
We say that a metric $d$ is {\it isotropic} if 
$$ \lim_{r\to 0} \inf_{|x-z|=r} d(x,z)= \lim_{r\to 0} \sup_{|x-z|=r} d(x,z)$$
for every $x$. 

\item[(ii)] 
The metric $d$ is called {\it locally star--shaped} if 
for every $x\in\X$ there exists an $r_0>0$ such that $B_d(x,r)$ is 
star--shaped with respect to the center of the ball, $x$, for every $r<r_0$. 
(A set $K$ is star--shaped with respect to $x$ if every ray emanating at 
$x$ intersects $\partial D$ exactly once.)

\item[(iii)]
The metric $d$ is called {\it locally convex} if for every $x\in\X$ there 
exists an $r_0>0$ such that $B_d(x,r)$ is convex for every $r<r_0$.
\end{itemize}
}\loppu

\alku\lem\label{isotropyLem} If $f_x(y):=M(x,y)$ is continuous at $x$ 
for every $x\in(0,\infty)$ then $\rho_M$ is isotropic.
\loppu

\proof Fix a point $x\in\X$. If $x=0$ then $\rho_M(x,z)=\rho_M(x,y)$ for 
every $|z|=|y|$. Let then $x\ineq 0$. If $f_x(x)=0$ then 
$\lim_{r\to 0} inf_{|x-z|=r} \rho_M(x,z) = \infty$ and $\rho_M$ is isotropic 
at $x$. Let then $c:= f_x(x)>0$. For every $0<\epsilon<f_x(x)/2$ 
there exists a neighborhood $U$ of $x$ such that $|f_x(y)-f_x(z)|\le \epsilon$. 
Then 
$$ \sup_{|x-z|=r} \rho_M (x,z)- inf_{|x-z|=r} \rho_M(x,z) \le 
 {r \over c-\epsilon} - {r \over c+\epsilon} \le 2r\epsilon /c$$
for every $0<r<d(\partial U, x)$ (here $d$ refers to the Euclidean distance). 
$\square$

\alku\rem{\rm 
It is possible that $\rho_M$ is isotropic even when $M$ is not continuous. For instance 
if $M(x,y)=xy$ for $x+y>0$ and $M(0,0)=1$ then $\rho_M$ is an isotropic metric, but 
clearly $M$ is not continuous at the origin. This example is due to Pentti J\"arvi.
}\loppu

\alku\lem\label{*shape} Let $\X$ be an inner product space. If $M$ is 
moderately increasing and $\rho_M$ is a metric then it is locally star--shaped.
\loppu

\proof Let us consider balls centered at $z$. Since the case $z=0$ is trivial,
we assume $z\ineq 0$. The case $M\ident 0$ is also trivial and then, since 
$M$ is moderately increasing, $M(x,y)>0$ for every $xy>0$. 

Let $r$ be a unit vector. Now if $\rho_M(z,z+sr) = s/M(z,z+sr)$
is increasing in $s>0$ for some range independent of the direction of $r$ 
then we are done. If $\Vert z+sr\Vert$ is decreasing in $s$, then $\rho_M$ 
is the product of two positive increasing factors, $s$ and $1/M(z,z+sr)$, 
and is hence itself increasing. 

If $\Vert z+sr\Vert$ is increasing in $s$, we write 
$$\rho_M(z,z+sr) = {s\over M(z,z+sr)} = {\Vert z+sr\Vert \over M(z,z+sr)}
{ s\over \sqrt{ \Vert z\Vert^2 +s^2 + 2s(r,z)} },$$
where $(r,z)$ denotes the inner product of $r$ and $z$. The first factor is 
increasing by the moderation part of the moderately increasing
condition of $M$. The second factor is increasing provided 
$\Vert z\Vert^2 \ge - s^2(r,z)$. Since $(r,z)\ge -\Vert z\Vert$ 
$\rho_M(z,z+sr)$ is increasing for $s\le \sqrt{\Vert z\Vert }$. 

Since $M$ is moderately 
increasing $M(x,y)$ is bounded from above in $B^n(z,s)$, say by $c_z$. 
Then $B_M(z,s/c_z)\subset B^n(z, s)$ and hence $B_M(z,s/c_z)$ is 
star--shaped. $\square$

\bigskip

The local star--shaped condition says that that the metric increases locally 
when we move away from the point (in the Euclidean metric), the isotropy condition
says that it does so equally fast in every direction. Both of these facts follow from 
the convexity result that we prove next.  

\alku\lem\label{locConv} Let $\X=\Rn$, $M$ be moderately 
increasing and $\rho_M$ be a metric. 
Assume also that $M(x,\cdot) \in C^2(\R^+)$. Then $\rho_M$ is locally convex.
\loppu

\proof Without loss of generality we may assume that $\X=\R^2$ since $B_M(z,r)$ is 
formed by rotating a two dimensional disk $B_M(z,r)\cap \R^2$ about the axis 
$tz$. Let us consider disks about $ze_1$, in particular, the locus of points $(x,y)$ with 
$\rho_M((x,y),z)=r>0$, i.e\abrv points for which the following equation holds:
\be\label{locConvE1} {\sqrt{(x-z)^2 + y^2} \over M(\sqrt{x^2+y^2},1)} = r. \ee

We will first show that if $y>0$ then $d^2y/dx^2 <0$. Let us denote 
$M(\sqrt{x^2+y^2},z)$ by $M$, $dM(w,z)/dw$ by $M'(w)$ and $d^2M(w,z)/dw^2$ by $M''(w)$. 
We multiply (\ref{locConvE1}) with $M$ and square it, then we differentiate 
with respect to $x$:
$$ yy' +x-z = r^2MM'(x+yy')(x^2+y^2)^{-1/2}.$$
From this it follows that
\be\label{locConvE2} yy'= \left( 1- {r^2MM'\over \sqrt{x^2+y^2}}\right)^{-1} -x.\ee
Differentiating again gives
$$ (y')^2 + yy'' = r^2 \left( 1- {r^2MM'\over \sqrt{x^2+y^2}}\right)^{-2}
\left( {(M')^2 +MM'' \over \sqrt{x^2+y^2}} - {MM' (x^2+y^2)^{-1}\over \sqrt{x^2+y^2} - 
r^2MM'}\right) -1.$$
By choosing $r$ sufficiently small, we may assume that $(x,y)\in B^2(1,\delta)$
for arbitrary given $\delta>0$. Then $\sqrt{x^2+y^2}\in [1-\delta,1+\delta]$
and there exists a constant $c $ such that $M(z),M'(z),M''(z) \le c$
for $z\in [1-\delta,1+\delta]$ since $M\in C^2(\R^+)$. If we choose $\delta\ge 1/2$
we also have $M(z)M'(z)/z\le 2c^2$ so that 
$$ yy'' \le 4r^2 c^2  (1- 2r^2 c^2)^{-2} -1 - (y')^2 .$$
Since $c$ is a constant it follows that $y''<0$ for sufficiently small $r>0$. 

Call the curve formed by the points which satisfy (\ref{locConvE1}) $\gamma$. 
If $\gamma$ could be parameterized by $(x,y(x))$ in Euclidean coordinates then the 
fact that $d^2y/dx^2 <0$ would imply that it is convex. 
Suppose that $\gamma$ can not be parameterized by $(x,y(x))$ (as shown in 
Figure 1). Then some half-line $K:= \{(x,t)\colon t\in\Rp\}$ intersects $\gamma$ at 
least twice. It follows that $dy/dx=\infty$ for some point $w$ in the upper half-plane. 
However, we see from (\ref{locConvE2}) that this is only possible for $y=0$, provided 
$r$ is small enough. Hence $w$ is in the $e_1$-axis, which is impossible. It follows 
that that $\gamma$ can be parameterized by $(x,y(x))$ and that the area under the 
curve is convex. Since $B_M(ze_1,r)$ is symmetric with respect to the $e_1$-axis it 
is convex, as well. $\square$

\begin{figure}
  \center{\epsfig{file=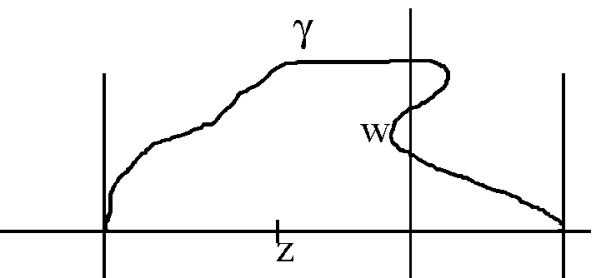,width=10cm}}
  \caption{Proof of Lemma \ref{locConv}} 
\end{figure}

%\alku\cor\label{pqLocConv} The $(p,q)$--relative metric
%is locally convex if and only if $p<\infty$. 
%\loppu

% replaced reference
\alku{\bf Proof of Theorem 1.2(2). \hskip 0.5truecm } {\rm
It is immediately clear that ${A_p(x,1)}^q \in C^2(\R^+)$ if $p<\infty$.

For $p=\infty$ we have $\rho_{\infty,q}(x,1)=\vert x-1\vert /\max\{1,\vert x \vert^q\}$.
Let us write $S(1,r)$ in polar coordinates about 1. Then $s(\theta)=r$ for 
$\cos\theta\ge r/2$ and 
$$ s^2 = r^2 (s^2+1-2s\cos\theta)^q $$
for $\cos \theta \le r/2$. It follows that for $\cos \theta < r/2$ we have
$$ s s'= r^2 q (s^2+1-2s\cos\theta)^{q-1} (ss'-s'\cos\theta +s\sin\theta).$$
Denote $\theta_0=\arccos(r/2)$. Since $s\to r$ as $\theta \to \theta_0^+$ 
($\theta$ approaches$\theta_0$ from above), we have 
$$ \lim_{\theta \to \theta_0^+}  s' = r^2q\sqrt{4-r^2}/(2-r^2q) >0.$$
Since $\lim_{\theta \to \theta_0^-}  s' =0$, the point $(2\cos\theta_0,\theta_0)$, 
will be an inner corner of $S(1,r)$ for every $r>0$, which means that 
$S(1,r)$ is not convex. $\square$
}\loppu

\alku\rem{\rm If a metric $d$ is locally star--shaped, isotropic or locally 
convex then so are $\log (1+d)$, $\arsh\, d$ and $\arch(1+d)$. Moreover, 
provided that $M$ is continuous these properties are also carried over to 
the generalized relative metrics considered in Section 6 of \cite{Ha}. 
}\loppu

%%%%%%%%%%%%%%%%%%%%%%%%%%%%%%%%%%%%%%%%%%%%%%%%%%%%%%%%%%%%%%%%%%%%%%%%%%%%%%%%%%%
%%%%%%%%%%%%%%%%%%%%%%%%%%%%%%%%%%%%%%%%%%%%%%%%%%%%%%%%%%%%%%%%%%%%%%%%%%%%%%%%%%%
%%%%%%%%%%%%%%%%%%%%%%%%%%%%%%%%%%%%%%%%%%%%%%%%%%%%%%%%%%%%%%%%%%%%%%%%%%%%%%%%%%%

\vskip 30pt 
\section{The hyperbolic metric and limitations of our approach} 

In this section, we will introduce the hyperbolic metric, show how our method 
can be used to generalize the hyperbolic metric in one setting but not in another. 
We use a separate method to deal with the latter case, thus solving a problem from
\cite[Remark 3.29]{Vu}. 

The hyperbolic metric can be defined in several different ways, for a fuller account the 
reader is referred to an introductory work on hyperbolic geometry, for instance \cite[Section 2]{Vu}. One possible definition of the hyperbolic metric, $\rho$, is
\be\label{hypMetrDef} \rho(x,y):= 2 \arsh \left({\vert x-y\vert \over 
\sqrt{1-\vert x\vert^2}\sqrt{1-\vert y\vert^2}}\right)\ee
for $x,y\in B^n$. An important property of the hyperbolic metric 
is that it is invariant under M\"obius mappings of $B^n$. The groups formed by these
M\"obius mappings is denoted by $\GM$. 

\alku\lem Let $M(x,y)=f(x)f(y)$ with $f(0)=1$ be such that $\rho_M$ is a metric. 
Then $\rho_M$ is invariant under all mappings in $\GM$ if and only if 
$f(x)=\sqrt{1-x^2}$.
\loppu

\alku\rem{\rm Note that here $f(x)$ is defined only for $x\in[0,1)$. Therefore $\rho_M$ 
is not exactly an $M$--relative metric in the sense defined in Section 1. The 
interpretation is nevertheless clear; strictly speaking we could extend $f$ by setting 
$f(x)=0$ when it was not previously defined and relying on the conventions regarding 
$\infty$.
}\loppu

\proof The "if" part say essentially that the hyperbolic metric is M\"obius invariant, as is seen from (\ref{hypMetrDef}). and is hence clear, see e.g\abrv \cite[2.49]{Vu}. 
Assume, conversely, that $\rho_M$ is invariant under all mappings in $\GM$. 

Fix $0<r<1$ and set $d:= r\sqrt{1-r^2}$. Then $d<2r$ and we may choose points 
$x,y\in B^n$ with $\vert x\vert =\vert y\vert =r$ and $\vert x-y\vert =d$. Let 
$g$ be a M\"obius mapping in $\GM$ which maps $y$ onto the origin. It 
follows from \cite[2.47]{Vu}, that $\vert g(x)\vert=r$. Hence by M\"obius invariance,
$$ {d \over f(r)^2} = {\vert x-y\vert \over f(\vert x\vert)^2} =
{\vert g(x)-0\vert \over f(\vert g(x)\vert) f(0)} = {r\over f(r)}$$
hence $f(r)=d/r=\sqrt{1-r^2}$. $\square$

\bigskip

The classical definition of the hyperbolic metric makes sense only in the unit ball and
domains M\"obius equivalent to it (for $n\ge 3$). There are however various generalizations of the hyperbolic metric to other domains. The best known of these 
is probably the quasihyperbolic metric that we met in Section 4. 
The quasihyperbolic metric is within a factor of 2 from 
the hyperbolic metric in the domain $B^n$ (\cite[Remark 3.3]{Vu}). 

Seittenranta's cross ratio metric is another
generalization of the hyperbolic metric, with the advantage, that it equals the 
hyperbolic metric in $B^n$. The reader may recall that we showed in \cite{Ha}, 
Corollary 6.5, that Seittenranta's metric can be interpreted as $\delta_G^{-\infty}$ 
in the one--parameter family $\delta_G^p$, 
$$ \delta_G^p(x,y) := \log \{1 + \rho'_{M,G}(x,y)\}$$
with $M= \max\{1, 2^{-1/p}\} A_p$, where $A_p$ is the power-mean, 
$$ A_p(x,y):=\left( {x^p+y^p \over 2} \right)^{1/p}$$
for $p\in (-\infty,0)\cup (0,\infty)$ and 
$$ A_{-\infty}(x,y)=\min\{x,y\},\ A_0(x,y):= \sqrt{xy}
{\rm\ and\ } A_{\infty}(x,y)=\max\{x,y\} $$ 
defined for $x,y\in\Rp$. Here 
$$ \rho'_{M,G}(x,y) = \sup_{a,b\in \partial G} {1 \over M( \vert x,y,a,b\vert, \vert x,y,b,a\vert) },$$
where 
$$ \vert a,b,c,d\vert :={q(a,c)q(b,d)\over q(a,b) q(c,d)}$$
denotes the cross--ratio of the points $a,b,c,d\in {\overline \Rn}$.

Seittenranta's metric is the generalization of the logarithmic expression for 
the hyperbolic metric given in \cite[Lemma 8.39]{Vu}. We now move on to study a 
generalization starting from the expression based the hyperbolic cosine 
(\cite[Lemma 3.26]{Vu}):
\be\label{rhoDef}  \rho_G(x,y):= \arch\{ 1 + \sup_{ a,b \in \partial G } 
\vert a,x,b,y \vert \vert a,y,b,x \vert/2 \}.\ee
This can be expressed as 
$$ \rho_G(x,y):= \arch\{ 1 + (\rho'_{A_0,G}(x,y))^2/2\},$$
with $A_0(x,y):=\sqrt{xy}$. 

We note that by \cite[Corollary 6.5]{Ha} we know that
$$ \log\{ 1 + \rho'_{A_0,G}(x,y)\}$$
is a metric provided $\card \partial G \ge 2$. Hence by \cite[Remark 3.7]{Ha} we already 
know that 
$$ \arch\{ 1 + \rho'_{A_0,G}(x,y)\}$$
is a metric when $\card \partial G \ge 2$. Hence one might speculate that the area 
hyperbolic cosine representation of the hyperbolic metric could be generalized to the 
one--parameter family 
$$ \rho_G^p(x,y):= \arch\{ 1 + (\rho'_{A_0,G}(x,y))^p/p\}.$$
In what follows we will however restrict our attention to the case $p=2$. 

Since this quantity has previously attracted some interest, we state some of its basic
properties and give an independent proof that it is in fact metric in most domains:

\alku\label{basicsThm}\thm {\rm (\cite[3.25 \& 3.26]{Vu})} 
\begin{itemize}
\item[(i)] $\rho_G$ is M\"obius invariant.
\item[(ii)]$\rho_G$ is monotone in $G$, that is, if $G\subset G'$ then 
$\rho_{G'}(x,y)\le \rho_G(x,y)$ for all $x,y\in G$. 
\item[(iii)] $\rho_G(x,y)\ge \cosh\{(q(\partial G)q(x,y))^2\} -1 $. 
\item[(iv)] For $G=B^n$ and $G=H^{n+}$ (the upper half-plane), $\rho_G$ equals the hyperbolic metric. 
\end{itemize}
\loppu

Note that $\rho_G$ is almost a generalized relative metric, indeed, we have 
$$ \rho_{\Rn\setminus\{0\}}(x,y):= \arch\left( 1 + {|x-y|^2\over 2|x\Vert y|}\right).$$
(Note that here $\Rn\setminus\{0\}$ has the boundary points $0$ and 
$\infty$ in $\Rnbar$.)
This expression differs from a generalized relative metric (essentially) 
only by the exponent 2 of $\vert x-y\vert$. However, because of this difference 
the question of whether it is a metric does not lend itself to the 
generalized metric approach of Section 6, \cite{Ha}. 

\alku\label{Thm1}\thm The quantity $\rho_G$ defined in (\ref{rhoDef}) is a metric
for every open $G\subset {\overline \Rn}$ with $\card\, \partial G\ge 2$.
\loppu

\proof It is clear that $\rho_G$ is symmetric in its arguments. That $(x,x)$
are the only zeros of $\rho_G$ is also evident. Moreover, as $\card  \partial G\ge 2$,
$\rho_G$ is finite. It remains to check that it satisfies the triangle inequality.

Since the supremum in the definition (\ref{rhoDef}) is over a compact set 
(in $\Rnbar$) it is actually a maximum. Fix $x$, $y$ and $z$ in $G$. 
Let $a,b\in \partial G$ be points such that
$$ \ch \rho_G(x,y)= 1 + \vert a,x,b,y \vert \vert a,y,b,x \vert/2.$$
Define $s(a,x,y,b):= \vert a,x,b,y \vert \vert a,y,b,x \vert/2$.
Now 
$$ \arch(1+s(a,x,z,b)) \le \rho_G(x,z),\ \arch(1+s(a,z,y,b)) \le \rho_G(z,y).$$
Hence it suffices to prove
\be\label{1} \arch(1+s(a,x,y,b)) \le \arch(1+s(a,x,z,b)) + \arch(1+s(a,z,y,b)).\ee

Since $s$ is conformally invariant, we may assume that $a=0$ and $b=\infty$.
Denote 
$$ s:= s(0,x,z,\infty)/2,\ t:= s(0,z,y,\infty) /2,\ u:=s(0,x,y,\infty) /2.$$
It follows that
\be\label{orgstu} s={\vert x-z\vert^2 \over 2\vert x\vert \vert z \vert },
 \ t={\vert z-y\vert^2 \over 2\vert z\vert \vert y \vert },
 \ u={\vert x-y\vert^2 \over 2\vert x\vert \vert y \vert }.\ee
For fixed $x$ and $y$ it is clear that we can move the point $z$ so that both
$s$ and $t$ get smaller if $\vert z\vert\le \min\{\vert x\vert, \vert y\vert\}$
(since $s= (x/z)+(z/x) - 2\cos\theta$ is increasing in $z$ for $z\le x$, and 
similarly for $t$). 
Hence we may assume that $\vert z\vert\ge \min\{\vert x\vert, \vert y\vert\}$. 
Similarly, if $\vert z\vert > \max\{\vert x\vert, \vert y\vert\}$ we can 
decrease $s$, $t$ for fixed $x$ and $y$, hence 
we may also assume that $\vert z\vert\le \max\{\vert x\vert, \vert y\vert\}$. 
Unless $\widehat{x0y}=\pi$ we may also assume that $z$ lies within this angle. 
Otherwise we may apply the transformations shown in Figure 2 (keeping $x$, $y$ and 
$|z|$ fixed and rotating or mirroring $z$ according to where it started.)

\begin{figure}
  \center{\epsfig{file=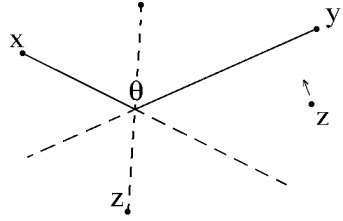,width=4.75cm}}
  \caption{The point $z$ is between $x$ and $y$}
\end{figure}

Since $\ch$ is increasing, we apply it to both sides of (\ref{1}) and use
$$ \ch(a+b)=\ch(a)\ch(b) + \sh(a)\sh(b)$$
to conclude that (\ref{1}) is equivalent to
\be\label{rhoTrigE0} u \le s+t+st+\sqrt{s^2+2s}\sqrt{t^2+2t}.\ee
Getting rid of the square--root, this equation is implied by
$$ s^2+t^2+u^2 \le 2(st+su+tu+stu) $$
which is equivalent to
\be\label{2} (u-s-t)^2 \le (4+2u)st.\ee

Let us assume without loss of generality that $z=1$. Assume, for the time being, that 0, $x$, $y$ and $1$ are co--linear and that $x>1>y>0$. Then
\be\label{3} s={1\over 2}\left(\sqrt{x}-{1\over \sqrt{x}}\right)^2,\ 
             t={1\over 2}\left(\sqrt{y}-{1\over \sqrt{y}}\right)^2,\
             u={1\over 2}\left(\sqrt{y\over x}- \sqrt{x\over y}\right)^2.\ee
Inserting these into (\ref{2}) gives
$$ \vert x+y+{1\over x}+ {1\over y}-{y\over x}-{x\over y}-2\vert \le
\left(\sqrt{y\over x}+\sqrt{x\over y}\right) \left(\sqrt{x}-\sqrt{1\over x}\right) 
\left(\sqrt{1\over y}- \sqrt{y} \right),$$
which is actually an equality.

Let us now consider the general case in which 0, $x$, $y$ and 1 are no longer necessarily
co--linear. Denote $s$, $t$ and $u$ from (\ref{3}) by $s_0$, $t_0$ and $u_0$, 
respectively and let $s$, $t$ and $u$ be as in (\ref{orgstu}). Denote
$$ \delta_s:=s-s_0=(1-\cos \theta),\ \delta_t:=t-t_0=(1-\cos \phi),\ 
\delta_u:=u-u_0=(1-\cos (\theta+\phi)),$$
where $\theta := \widehat{x01}$ and $\phi:= \widehat{10y}$.
Inserting $s=s_0 + \delta_s$ etc\abrv into (\ref{2}) and canceling the equality 
$(s_0+t_0-u_0)^2=2(2+u_0)s_0t_0$ leads to
$$ 2(s_0+t_0-u_0)(\delta_s+\delta_t-\delta_u) + 
(\delta_s+\delta_t-\delta_u)^2\le 2\delta_u st + 2(2+u_0)(t_0 \delta_s+ 
s_0 \delta_t + \delta_s \delta_t)$$
which is equivalent to
\be\label{rhoTrigE1} (2s_0+\delta_s)(\delta_s-\delta_t-\delta_u) + 
(2t_0+\delta_t) (\delta_t-\delta_s-\delta_u) + 
(2u_0+\delta_u) (\delta_u-\delta_s-\delta_t) \le $$
$$\le 2(stu-s_0t_0u_0).\ee

We will first show that 
\be\label{rhoTrigE2} 
\delta_s(\delta_s-\delta_t-\delta_u) + \delta_t (\delta_t-\delta_s-\delta_u) + 
\delta_u (\delta_u-\delta_s-\delta_t) \le 0. \ee
Note first that $\delta_s\ge 0$, $\delta_t\ge 0$ and $\delta_u\ge 0$. Now either 
all the parenthesis are negative or $\delta_u-\delta_s-\delta_t\ge 0$, since 
$\delta_u\ge \delta_s, \delta_t$. In the latter case the left hand side of the 
inequality is increasing in $\delta_u$. Since $\delta_s$, $\delta_t$ and $\delta_u$ 
are squares of the sides of a triangle we see that 
$$\delta_u \le \delta_s+ \delta_t + 2\sqrt{\delta_s\delta_t}.$$
Hence it suffices to check (\ref{rhoTrigE2}) for the maximal $\delta_u$, in which 
case it is an equality. 

Let us then continue from (\ref{rhoTrigE1}), using (\ref{rhoTrigE2}), rearranging and 
dividing by 2:
$$ \delta_s(s_0-t_0-u_0) + \delta_t(t_0-s_0-u_0) + \delta_u(u_0-s_0-t_0) \le 
stu-s_0t_0u_0.$$
Since $\delta_s,\delta_t\ge 0$ it follows that $stu-s_0t_0u_0 \ge s_0t_0\delta_u$. 
We will then complete the proof by showing that 
$$ \delta_s(s_0-t_0-u_0) + \delta_t(t_0-s_0-u_0) + \delta_u(u_0-s_0-t_0-s_0t_0)\le 0.$$
We may assume that (\ref{rhoTrigE0}) holds with equality, hence 
$$ u_0 = s_0+t_0+s_0t_0+ \sqrt{s_0^2+2s_0}\sqrt{t_0^2+2t_0}.$$
Then it suffices to show that 
\be\label{rhoTrigE3} (\delta_u-\delta_s-\delta_t) \sqrt{s_0^2+2s_0}\sqrt{t_0^2+2t_0} \le 
2(t_0\delta_s + s_0\delta_t) + (\delta_s+\delta_t)s_0t_0.\ee
By the formula for the cosine of a sum we have, from the definition,
$$ \delta_u= \delta_s +\delta_t -\delta_s\delta_t +\sqrt{2\delta_s-\delta_s^2}\sqrt{2\delta_t-\delta_t^2} \ge 
\delta_s +\delta_t  + \sqrt{2\delta_s-\delta_s^2}\sqrt{2\delta_t-\delta_t^2}.$$
Then (\ref{rhoTrigE3}) follows if we can show that 
$$ \sqrt{2\delta_s-\delta_s^2}\sqrt{2\delta_t-\delta_t^2} \sqrt{s_0^2+2s_0}\sqrt{t_0^2+2t_0} \le 2(t_0\delta_s + s_0\delta_t) + (\delta_s+\delta_t)s_0t_0.$$ 
Let us square this equation and subtract $2\delta_s\delta_t st(2+s)(2+t)$ from 
both sides:
$$ (2-2(\delta_s+\delta_t) + \delta_s\delta_t) \delta_s \delta_t (2+s_0)(2+t_0)s_0t_0 
\le \delta_s^2t_0^2(2+s_0)^2 + \delta_t^2 s_0^2 (2+t_0)^2.$$
Divide both sides by $\delta_s \delta_t (2+s_0)(2+t_0)s_0t_0$:
$$ 2-2(\delta_s+\delta_t) + \delta_s\delta_t \le a+1/a,$$
where 
$$ a:= {\delta_s (2+s_0) t_0 \over \delta_t s_0 (2+t_0)} $$ 
(this is OK, since the cases where $\delta_t=0$ or $s_0=0$ are trivial.) 
Now then $a+1/a\ge 2$ (by the arithmetic-geometric inequality, for instance) 
so it suffices to show that $\delta_s\delta_t \le 2(\delta_s+\delta_t)$ or
equivalently,
$$ {1\over 2} \le {1\over \delta_s}+ {1\over \delta_t}.$$ 
But since $\delta_s,\delta_t\le 2$ directly from the definition, this is 
clear. $\square$

\alku\ack{\rm I would like to thank and Matti Vuorinen 
for numerous comments and suggestions as well as Glen D. Anderson and 
Pentti J\"arvi for their comments on earlier versions of this manuscript.
}\loppu

%%%%%%%%%%%%%%%%%%%%%%%%%%%%%%%%%%%%%%%%%%%%%%%%%%%%%%%%%%%%%%%%%%%%%%%%%%%%%%%%%%
%%%%%%%%%%%%%%%%%%%%%%%%%%%%%%%%%%%%%%%%%%%%%%%%%%%%%%%%%%%%%%%%%%%%%%%%%%%%%%%%%%
%%%%%%%%%%%%%%%%%%%%%%%%%%%%%%%%%%%%%%%%%%%%%%%%%%%%%%%%%%%%%%%%%%%%%%%%%%%%%%%%%%

\end{document}